%% file: banach.tex
\catcode`\@=11
\documentclass[11pt]{amsart}

\title[On flat trigonometric sums and simple Lebesgue spectrum]
      {On flat trigonometric sums and ergodic flow with simple~Lebesgue~spectrum}
\author{A.\,A.\,Prikhod'ko}
\date{11.11.2009}

\usepackage{amsmath}
\usepackage{amscd}
\usepackage{amsfonts}
\usepackage{amssymb}

\usepackage[dvips]{graphicx}

\usepackage{sasha_prih}
\usepackage{rank1}

\usepackage{amsthm}  

\newtheorem{thm}{Theorem}[section]
\newtheorem{lem}[thm]{Lemma}

\newtheorem{cor}{Corollary}

\theoremstyle{definition}
\newtheorem{defn}{Definition}[section]

\theoremstyle{remark}
\newtheorem{rem}{Remark}

\begin{document}

\begin{abstract}
A complex polynomial ${P(z) = c_0 + c_1 z +\ldots+ c_n z^n}$  
is called {\it unimodular\/} if ${|c_j| = 1}$, ${j = 0,\ldots,n}$. 
Littlewood 
asked the question (1966) on how close a unimodular polynomial
come to satisfying ${|P(z)| \approx \sqrt{n+1}}$ if ${n \ge 1}$\,? 
In this paper 
we show that for a given ${0 < a < b}$ and ${\eps > 0}$ 
there exist trigonometric sums 
${\cP(t) = n^{-1/2} \sum_{j=0}^{n-1} \exp(2\pi i\, t\omega(j))}$ with 
a real frequency function $\omega(j)$
which are $\eps$-{\it flat\/} on segment $[a,b]$ acording to the norm in $L^1([a,b])$ 
(as well as in $L^2([a,b])$). 
We apply this method to construct a dynamical system 
having simple spectrum and Lebesgue spectral type in the class of rank-one flows. 

The work is supported by grant \No\,NSh-3038.2008.1.
\end{abstract}

\maketitle

%
\section{Introduction}

In this paper we study the spectral properties of a class of ergodic dynamical systems satisfying
certain approximation properties and called {\it system of finite rank}. 
In brief we say that a measure preserving transformation $T$ on a Lebesgue probability space 
has {\it rank one\/} if there exists a sequence of Rokhlin towers 
${TB \sqcup \ldots \sqcup T^{h_n}B}$ approximating any measurable set. 
Rank one dynamical systems are known to have simple spectrum. 
We show that there exist rank one flows of measure preserving transformations 
having Lebesgue spectral type which gives positive answer to Banach question 
on simple Lebesgue spectrum for ergodic $\Set{R}$-actions. 
The method we use is actually based on purely analytic problem on 
existence of flat polynomials and trigonometric sums which goes back to 
classical Littlewood question about how close a unimodular polynomial
$$
	P(z) = c_0 + c_1 z +\ldots+ c_n z^n, \qquad |c_j| = 1, 
$$ 
come to satisfying ${|P(z)| \approx \sqrt{n+1}}$ if ${n \ge 1}$\,? 
In particular case of restriction ${c_j = \pm 1}$ the question still remains open. 
Bourgain \cite{Bourgain} discovered 
that the spectral type of a rank one transformation 
is given by generalized Riesz product, and 
as an ingredient it contains polynomials with coefficients $0$ and~$1$ 
which are similar and connected to polynomials with Littlewood type coefficient constraints. 
The purpose of our investigation is 
to construct $\eps$-flat trigonometric sums in order to 
find a class of rank~one flows having simple Lebesgue spectrum. 

\subsection{Spectral invariants of dynamical systems}

Given $T \Maps X \to X$ a~measure preserving invertible transformation 
on a standard Lebesgue space %
$(X,\cX,\mu)$ we consider a unitary operator $\Hat T$ on $L^2(X,\cX,\mu)$ 
called {\it Coopman operator\/} acting as shift along the trajectory, 
${\Hat T f(x) = f(Tx)}$. Since constants are invariant under $\Hat T$ the operator $\Hat T$ 
is usually restricted to the space $L^2_0(X)$ of functions with zero mean, ${\int f \,d\mu = 0}$. 
For each ${f \in L^2_0(X)}$ consider the {\it spectral measure\/} $\sigma_f$ on the unit circle 
${S^1 = \{z \in \Set{C} \where |z| = 1\}}$ defined by 
\begin{equation}
	\int_{S^1} z^k\,d\sigma_f = \scpr<T^kf,f>. 
\end{equation}
According to spectral theorem $\Hat T$ is uniquely determined up to spectral equivalence 
by two invariants, the spectral type $\sigma$ on the unit circle $S^1$ and 
the multiplicity function $\Bar m(z)$.  
For any $f$ measure $\sigma_f$ is absolutely 
continuous with respect to the {\it spectral type measure\/} $\sigma$, 
we denote it as ${\sigma_f \ll \sigma}$, and $\sigma$ is the minimal meaasure 
having this property up to the following measure equivalence: 
${\la_1 \sim \la_2}$ iff ${\la_1 \ll \la_2}$ and ${\la_2 \ll \la_1}$ (see \cite{Halmos}). 

We say that $T$ has {\it simple spectrum\/} if 
$T$ possesses a {\it cyclic vector}, an element ${h \in L^2_0(X)}$ such that 
${Z(h) = L^2_0(X)}$, where ${Z(h) = \CSpan\{T^kh \where k \in \Set{Z}\}}$. 
The equivalent way to define simple spectrum property is to require 
${\Bar m(z) = 1}$ for the multiplicity function. 
Well known examples of simple spectrum maps are ergodic dynamical system with 
purely discrete spectrum studied by von Neuman (see \cite{Halmos}). In this case 
the spectral measure is supported on a discrete subgroup $\Lambda$ of $S^1$. 

A variety of constructions of dynamical systems with finite spectral multiplicity 
and different dynamic effects comes from {\it approximation theory\/} inspired by 
Rokhlin ideas on tower approximation and developed by Katok, Oseledec and Stepin 
(see \cite{KSt,Oseledec1,Stepin86}). 
%
%
It was discovered that spectral multiplicity of $r$-interval exchange maps is bounded by ${r-1}$, 
and further Robinson \cite{Robinson} proved that there exist interval exchange maps 
realizing maximal multiplicity of ${r-1}$. 
Ageev in a series of works (e.g.\ see \cite{Ageev88,Ageev05}) 
provided a set of constructions to get different multiplicity essential value set. 
Though the question on existence of an ergodic transformation 
having homogeneous spectral multiplicity ${m > 1}$ known as Rokhlin problem 
was open by 1999 when Ryzhikov and Ageev constructed first examples 
of ergodic maps with homogeneous spectral multiplicity two
(see \cite{AnosovOKratnostyah}).

\subsection{Rank one transformations and flows}

\begin{defn}
Let $T$ be a measure preserving invertible transformation on a Lebesgue space $(X,\cX,\mu)$. 
Given a set ${B \in \cX}$ and integer $h$ assume that 
$B,\,TB,\,T^2B,\,\ldots\,T^{h-1}B$ are disjoint. 
The set ${\cT_{h,B} = \bigcup_{j=0}^{h-1}T^kB}$ is called {\it Rokhlin tower\/} 
of height $h$, and $B$ is called the {\it base\/} of the tower. 
We can draw a tower as a sequence of {\it levels\/} $T^kB$ such that $T$ lifts 
$k$-th level $T^kB$ to the next level $T^{k+1}B$ (except the top level). 
We refer to the following partition as {\it level partition\/} of the tower 
\begin{equation}
	\xi_{h,B} = \{B,TB,\ldots,T^{h-1}B,X\sms \cT_{h,B}\}. 
\end{equation}
\end{defn}

The well-known Rokhlin--Halmos' lemma (see \cite{Halmos}) states that 
given aperiodic transformation $T$ for any height $h$ and arbitrary ${\eps > 0}$ 
there exists a Rokhlin tower $\cT_{h,B}$ satisfying ${\mu(\cT_{h,B}) > 1-\eps}$. 

\begin{defn}
A measure preserving invertible map $T$ is called {\it rank one transformation\/} if 
${\mu(\cT_{h_n,B_n}) \to 1}$
and there exists a sequence of Rokhlin towers $\cT_{h_n,B_n}$ approximating 
$\sigma$-algebra in the following sense. For any measurable set $A$ we can find 
$\xi_{h_n,B_n}$-measurable sets $A_n$ such that 
${\mu(A_n\syms A) \to 0}$ as ${n \to \infty}$. 
\end{defn}

Rank one transformations was introduced by Chacon \cite{Chacon}, Ornstein \cite{O}
and Fridman \cite{Friedman1}. Some classical ergodic maps, for example 
discrete spectrum maps appear to be rank one (see survey \cite{Ferenczi1}). 
A series of dynamical system with finite rank are found in the class 
of adic and substitution systems~\cite{VershikAddic}. 
%
%
In paper \cite{O} Ornstein constructed a randomized family of rank one maps 
known to be almost surely mixing which means that 
${\mu(T^{j}A \cap B) \to \mu(A)\,\mu(B)}$ as ${j \to \infty}$ 
for all measurable sets $A$ and $B$. 
For the class of Ornstein transformations Bourgain proved that 
almost surely they have purely singular spectrum \cite{Bourgain}. 
Today all known examples of rank one maps and flows have specatral type 
which is combination ${\sigma = \sigma_d + \sigma_s}$ of discrete $\sigma_d$ 
and singular $\sigma_s$ measures (see \cite{ChoksiNadkarni,Klemes,KlemesReinhold}). 

The question on existence of a map (or a flow) having simple spectrum 
(in orthocomplement to constants) 
and Lebesgue spectral type is known as Banach problem on simple Lebesgue spectrum 
(see \cite{KSF,Nadkarni}). 
Mathew and Nadkarni \cite{MathewNadkarni} 
answering a question by Helson and Parry \cite{HelsonParry} 
has constructed a transformation having a Lebesgue component in spectrum 
of multiplicity $2$ with a discrete component. 
Guenais \cite{Guenais} has found a criterion for existence of Lebesgue component 
in spectrum of generalized Morse dynamical systems. 

In this paper we focus on spectral type of rank one flows 
(see \cite{RyzhRankOneFlows,Prikhodko}). 

\begin{defn}
We say that a set $\cT_{h,B} = \bigcup_{t \in [0,h]} T^tB$ is 
a {\it Rokhlin tower\/} for a flow $T^t$ on a standard Lebesgue space $(X,\cX,\mu)$ 
if the sets $T^tB$ are disjoint and 
for any Borel set $J$ the union $\bigcup_{t \in J} T^tB$ is measurable. 
Let us define measurable partitions $\xi_{h,B}$ of the phase space $X$ 
into levels $T^tB$ of the tower and the set complementary to the tower.
\end{defn}

\begin{defn}
A flow $T^t$ is called {\it rank one\/} flow if there exist a sequence 
of Rokhlin towers $\cT_{h_n,B_n}$ such that ${\mu(\cT_{h_n,B_n}) \to 1}$ and 
for any measurable $A$ we can find $\xi_{h_n,B_n}$-measurable sets 
$A_n$ such that ${\mu(A_n \syms A) \to 0}$ as ${n \to \infty}$.
\end{defn}

We will assume that additionally partition $\xi_{h_{n+1},B_{n+1}}$ 
refines partition $\xi_{h_n,B_n}$. In this case we can define the rank one flow
by the following construction. 

\begin{defn}
{\it Cutting-and-stacking construction\/}: 
Consider a sequence ${q_n \in \Set{N}}$ and spacer parameters ${s_{n,k} \ge 0}$. 
At $n$-th step we cut the Rokhlin tower $\cT_{h_n,B_n}$ into $q_n$ 
equal subtowers called {\it columns}, add spacer $s_{n,k}$ to $k$-th columns 
and stack it together ordered from the first to the $q_n$-th. 
Repeating this procedure infinitely many times we construct a Lebesgue 
probability space if the folowing condition holds:
\begin{equation}\label{EqFiniteM}
	\prod_{n=1}^\infty \frac{h_{n+1}}{q_nh_n} < \infty, 
\end{equation}
where $h_n$ is the height of the $n$-th tower. 
A~point ${a \in X}$ moves up with unit velocity 
along any tower under action of the flow $T^t$ 
and it makes jump reaching the roof of the tower. 
\end{defn}

Equivalently cutting-and-stacking construction can be represented as follows. 

\begin{defn}
Consider segments ${X_n = [0,h_n]}$ and define a space $X$ 
to be the inverse limit of $X_n$ up to projections 
${\phi_n \Maps X_{n+1} \to X_n}$, 
\begin{equation}
	\phi_n(a_{n,k}+t) = t, \qquad 0 \le t < h_n, 
\end{equation}
where $a_{n,k} = \sum_{j < k}(h_n+s_{n,j})$, and ${\phi_n(u) = 0}$ otherwise. Namely,
\begin{equation}
	X_1 \stackrel{\phi_1}{\leftarrow} 
	X_2 \stackrel{\phi_2}{\leftarrow} \ldots 
	\stackrel{\phi_{n-1}}{\leftarrow} X_n \stackrel{\phi_n}{\leftarrow} \ldots X. 
\end{equation} 
If \eqref{EqFiniteM} holds we can normalize Lebesgue measure on $X_n$ 
so that $\phi_n$ become measure preserving and $X$ becomes Lebesgue probability space. 
A~point ${x \in X}$ is now a sequence 
\begin{equation}
	x = (x_1,x_2,\ldots) \quad \text{such that} \quad \phi_n(x_{n+1}) = x_n. 
\end{equation}
To define $T^t$ for a point $x$ we can write ${(T^tx)_n = x_n+t}$ if 
the point ${x_n+t}$ do not go out the right boundary of the segment $[0,h_n]$, 
elsewise we have to watch the point at some upper level ${n' > n}$. 
It can be easily seen that almost surely $x_n$ is $t$-far from the right boundary 
of $[0,h_n]$ starting from some index $n_0$, and the definition of the flow is correct. 
\end{defn}

\subsection{Unimodular polynomyals and trigonometric sums}

Complex polynomial 
\begin{equation}
	P(z) = c_0 + c_1 z +\ldots+ c_n z^n, \quad c_j \in \Set{C}, \quad 
	z \in \Set{C}, \quad |z| = 1, 
\end{equation}
is called {\it unimodular\/} if all the coefficients ${|c_j| = 1}$, 
${j = 0,1,\ldots,n}$. Let $\cK_n$ be the class containing 
unimodular polynomials of degree~$n$, and let $\cL_n$ be the class 
of polynomyals with coefficients $\pm 1$, 
\begin{equation}
	\cL_n = \left\{ Q(z) = \sum_{j=0}^{n} c_j z^j, \ c_j \in \{-1,+1\} \right\}.
\end{equation}
Consider a function $|P(z)|$ on a unit circle ${S^1 = \{|z|=1\}}$.
J.\,Littlewood (\cite{Littlewood}, cf.\ 
	\cite{ErdelyiLittlewoodType02,ErdelyiLTP99,Kahane80,BourgainBom})  
has formulated a question on how close a unimodular polynomial 
come to satisfying ${|P(z)| \approx \sqrt{n+1}}$ 
when $P$ range over the class $\cK_n$ or the class $\cL_n\,$? 
We say that a polynomial $P(z)$ is $\eps_n$-{\it flat\/} if 
\begin{equation}
	\left\|\,\frac1{\sqrt{n+1}}|P(z)|-1\,\right\| \le \eps_n
\end{equation} 
according to some norm. 
T.\,K\"orner \cite{Korner,Byrnes} 
has discovered a unimodular $P(z)$ with the property 
${ A \le (n+1)^{-1/2}|P(z)| \le B }$ for some absolute constants ${0 < A < 1 <  B}$. 
J.\,Kahane \cite{Kahane80} has constructed a pollynomial 
satisfying the following estimate at any point ${z \in S^1}$ 
\begin{equation}
	\bigl|\,(n+1)^{-1/2}|P(z)|-1\,\bigr| \le \eps_n, \qquad \eps_n = O(n^{1/17} \log n). 
\end{equation}
It is also an open question if there eexist a flat polynomial
in the class 
\begin{equation}
	\cM_n = \left\{ 
		P(z) = q^{-1/2}(z^{a_0}+z^{a_2}+\dots+z^{a_{q-1}}),\ 
		a_k \in \Set{Z},\ 
		a_k < a_{k+1} 
	\right\} 
\end{equation}
where ${n \ge 2}$.

\subsection{Flat trigonometric sums}

We will study trigonometric sums 
\begin{equation}
	\cP(t) = \frac1{\sqrt{q}} \sum_{y=0}^{q-1} \exp(2\pi i\, t \omega(y)), \qquad 
	\omega(y) \in \Set{R}, \qquad t \in \Set{R},  
\end{equation} 
where $\omega(y)$ is refered to as {\it frequency function}. 
The main target of our inversigating is the following special class of functions~$\omega(y)$. 

\begin{defn}
We call {\it exponential staircase frequency function\/} 
the function ${\omega \Maps \Set{Z} \to \Set{R}}$ 
given by the following expression:
\begin{equation}
	\omega(y) = E_0 + A e^{\epsilon y}, \qquad E_0,A,\epsilon \in \Set{R}. 
\end{equation} 
Observe that $\omega(y)$ is just a solution of the differential equation 
\begin{equation}
	\omega'' = \epsilon \omega' 
\end{equation} 
if we consider it for a while as a function on $\Set{R}$. 
\end{defn}

\begin{thm}
For given $0 < a < b$, ${\eps > 0}$ and ${\delta > 0}$ 
there exists $m_0$ such that for any ${m \ge m_0}$ 
there exists an infinite sequence $q_j$ generating 
trigonometric sums with exponential staircase frequency function
\begin{equation}
	\omega(y,m) = m\frac{q_j}{\eps^2}e^{\eps y/q_j}
\end{equation}
which are $\delta$-flat in $L^1([a,b])$.
\end{thm}

\subsection{Exponential staircase flows. Main result}

In this paper we study spectral type of a class of rank one flows 
given the following definition.

\begin{figure}[th]
	\begin{center}
		\includegraphics[width=100mm]{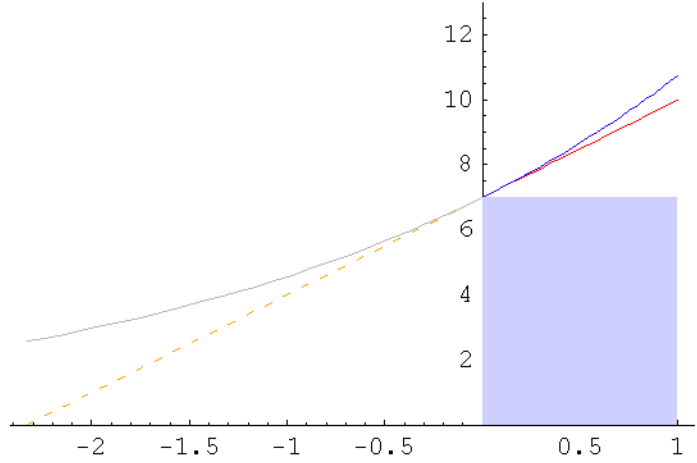}
	\end{center}
	\caption{Staircase and exponential staircase constructions: 
		Light blue rectangle is the $n$-th tower, red line is a roof by staircase construction,
		blue line is an exponential staircase roof.} 
	\label{pExpStaircase}
\end{figure}

\begin{defn}
We say that a rank one flow is given by {\it exponential staircase construction\/} 
if the roof function over the $n$-th tower is given by the graph of 
the discrete derivative ${\omega_n(y+1)-\omega_n(y)}$ of an exponential 
function $\omega_n(y)$, in other words if there exist $\eps_n$ and $m_n$ such that 
\begin{equation}
	h_n + s_{n,y} = \omega_n(y+1) - \omega_n(y), \qquad 
	\omega_n(y) = m_n\frac{q_n}{\eps_n^2}(e^{\eps_n y/q_n}-1), \qquad h_n = \frac{m_n}{\eps_n}. 
\end{equation}
Remark that $-1$ in the parenthesis is used to provide 
standard position of the segment $[0,h_n]$ in $\Set{R}$, namely, ${\omega_n(0) = 0}$, 
and if ${\eps_n \to 0}$ as ${n \to \infty}$ then 
\begin{equation}
	s_{n,y} \sim \omega_n'(y) - h_n = 
	m_n\frac{1}{\eps_n}e^{\eps_n y/q_n} - h_n = m_n\frac{1}{\eps_n}(e^{\eps_n y/q_n}-1). 
\end{equation}
It would be also interesting to observe that the spacer function $s_{n,y}$ 
ranges over ${[0,m_n(1+O(\eps_n)]}$. 
\end{defn} 

\begin{rem}
Terminology we use is explained by the note that spacer function $s_{n,k}$ 
approximates the spacer function ${\tilde s_{n,y} = \alpha y}$ which 
leads to {\it staircase construction\/} of a rank one flow, generalizing 
the concept of rank one staircase transformation introduced by 
Adams and Smorodinsky (see \cite{Adams}) providing the first explicit example 
of mixing rank one transformation (see also \cite{RyzhSpectrAndMixing}). 
\end{rem}

The main purpose of the paper is to study the spectral properties of 
exponential staircase constructions. It appears that in many cases 
the spectral type of such flows can be calculated in a certain sense, 
and we can find large variety of flows with Lebesgue spectral type 
having cardinality of continuum. 

\begin{thm}
Let us consider the exponential staircase construction of rank one flow $T^t$ 
given by the sequence of functions
\begin{equation}
	\omega_n(y) = m_n\frac{q_n}{\eps_n^2}e^{\eps_n y/q_n}, \qquad 
	h_n = \frac{m_n}{\eps_n}. 
\end{equation}
Suppose that for an $\alpha$ with ${0 < \alpha < 1/4}$ 
\begin{equation}
	h_n^{1+\alpha} \le q_n, \qquad m_n \le h_n^{1/2-\alpha}, \qquad 
	\sum_{n=1}^\infty \frac1{m_n^{1/4-\alpha}} < \infty, 
\end{equation}
and $q_n$ has the following property: 
\begin{equation}
	\forall k \in \Set{Z} \cap 
		\left[ \frac{m_n^{1-\alpha}}{\eps_n}, \frac{m_n^{1+\alpha}}{\eps_n} e^{\eps_n} \right] 
	\quad \left\| (q_n+1) \frac1{\eps_n}k\ln k \right\|_{\Set{Z}} \le \eps_n, 
\end{equation}
where
\begin{equation}
	\| x \|_{\Set{Z}} = \min\{|\ell - x\| \where \ell \in \Set{Z}\}. 
\end{equation}
Then the flow $T^t$ has Lebesgue spectral type. 
\end{thm}

{\it Sketch of the proof.} 
Consider a function $f$ which is measurable up to partition 
$\xi_{h_{n_0},B_{n_0}}$ into levels of the $n_0$-th tower 
in the cutting-and-stacking construction of rank one flow. 
We may consider $f$ as a function on $\Set{R}$ defining $f(t)$ 
to be the constant value $f$ keeps on the level $T^tB$. 
Using observation due to Bourgain on approximation of 
the spectral measure $\sigma_f$ we can say that it is given 
by a generalized Riesz product converging in weak sense, 
\begin{equation}
	\sigma_f = |\Hat f|^2 \cdot \prod_{n \ge n_0} |\cP_n(t)|^2, 
\end{equation}
where
\begin{equation}
	\cP_n(t) = \frac1{\sqrt{q_n}} \sum_{y=0}^{q_n-1} \exp 2\pi i\, t\omega_n(y), 
\end{equation}
and $\omega_n(y)$ is connected with the sequence of spacers as follows: 
${h_n + s_{n,y} = \omega_n(y+1)-\omega_n(y)}$. 
In~brief it easily follows from the observation that 
the lift $f_{(n)}$ of the function $f$ to the space ${X_{n} = [0,h_{n}]}$ 
and the lift $f_{(n+1)}$ to larger space $X_{n+1}$ are connected in the following simple way
\begin{equation}
	f_{(n+1)}(x) = \sum_{y=0}^{q_n-1} f_{(n)}(x-\omega_n(y)), 
\end{equation}
thus,
\begin{equation}
	\Hat f_{(n+1)}(t) = 
		\Hat f_{(n)}(t) \cdot \frac1{\sqrt{q_n}} \sum_{y=0}^{q_n-1} e^{2\pi i\, t\omega_n(y)}. 
\end{equation}

The idea is to find a sequence of trigonometric sums $\cP(t)$ 
which are $\delta_n$-flat on expanding intervals $(a_n,b_n)$, ${a_n \to 0}$, ${b_n \to \infty}$, 
with ${\delta_n \to 0}$ sufficiently fast, 
so that we can control convergence of the Riesz product 
to a regular density, and by ergodicity we see that 
$\sigma_f$ is absolutely continuous with respect to Lebesgue measure, 
and choosing appropriate $f$ we can assume that the density of $\sigma_f$ 
is not vanishing for any given intercal in $\Set{R}$.

\input st_phase.tex

\input proof_riesz.tex

The author is very grateful to D.\,V.\,Anosov, A.\,M.\,Vershik, 
A.\,A.\,Stepin, S.\,V.\,Konyagin, B.\,S.\,Kashin and M.\,Lemanczyk 
for discussions and remarks. I would like to thank 
V.\,V.\,Ryzhikov and J-P.\,Thouvenot for attention to this work and helpful discussions.
I would like to thank H.\ el Abdalaoui for fruitful discussions 
concerning generalized Riesz produsts and spectral properties of rank one flows.
The author is very grateful to P.\,Grigoriev for cooperation and help 
in some questions concerning harmonic analysis and polynomials, 
and to A.\,Egorov for helpful discussions and hints concerning 
analysis and number theory.

\input bibl09.tex
\end {document}

%% file: st_phase.tex
\section{Stationary phase method}

\subsection[Van der Corput's method]{Preliminary lemmas. Van der Corput's method} 
\label{sVanDerCorput}
  
Consider a smooth real-valued {\it frequency function\/} $f(x)$ 
on the real line $\Set{R}$ and look at the integral 
\begin{equation}
	\int_{x_0}^{x_1} e^{it f(x)}\,dx 
\end{equation}
An extremal point of the function $f$ defined by equation ${f' = 0}$ 
is called a {\it stationary phase\/} of~$f$. 
For~example, a frequency function ${f(x) = (a + cx^2)}$ has one stationary phase point 
${x_0 = -a/(2c)}$. 
Let us mention two well-known lemmas on oscillatory interals with quadratic $f$. 
  
\begin{lem}\label{Oi1} Oscillatory integral with quadratic frequency function 
and zero stationary phase can be estimated as follows 
\begin{equation}
	\int_{-1}^{1} e^{ikx^2}\,dx = 
		\sqrt{\frac{\pi}{k}} \exp \left( \frac{\pi i}{4} \right) + O\left(\frac1{k}\right), 
	\quad k \to \infty.
\end{equation}
\end{lem} 
  
\begin{lem}\label{Oi2} For real $a$, $c \not= 0$ and ${b > 0}$ 
\begin{equation}
	\int_{0}^{b} e^{it(a + cx^2)}\,dx = 
		A_0 \frac{e^{iat}}{2(|c|t)^{1/2}} - \frac{i}{2bct} e^{i(a+cb^2)t} + 
		  O \left( \frac1{b^3(ct)^2} \right)
\end{equation}
as $t \to \infty$, where $A_0$ is defined as 
\begin{equation}
A_0 = \int_0^\infty u^{-1/2} e^{iu\,\sgn(c)}\,du = e^{\frac14 \pi i\,\sgn(c)} \sqrt{\pi}. 
\end{equation}
\end{lem} 

%
Now let us consider a sum over an interval in the integer line 
\begin{equation}
	S = \sum_{a \le y \le b} e^{2\pi i f(y)}, \qquad y \in \Set{Z}, 
\end{equation}
where $f \in C^2([a,b])$. 
The following procedure is called {\it van der Corput's method}. 
The sum $S$ can be approximated by a special sum over stationary phase set 
\cite{Dieudonne,Titchmarsh,Liu}, namely, 
let $f'$ be non-decreasing on $[a,b]$, ${\alpha = f'(a)}$, ${\beta = f'(b)}$, then
\begin{equation}\label{vanDerCorputS}
	S = \sum_{\alpha < k < \beta} |f''(y_k)|^{-1/2} e^{2\pi i\, (f(y_k) - ky_k + 1/8)} + \cE, 
\end{equation}
where $y_k$~are solutions of the equation ${f'(y_k) = k}$, and $\cE$ is error term. 
  
\begin{lem} Suppose $\la_2 \le f''(y) \le \const \cdot \la_2$, 
${ |f^{(3)}(y)| \le \la_3 }$, ${\la_2,\la_3 > 0}$, 
and $\cE$ is defined in \eqref{vanDerCorputS} 
then (van der Corput, see\ \cite{Liu} and \cite{Titchmarsh} Theorem 4.9) 
\begin{equation}
\cE = O(\la_2^{-1/2}) + O(\ln(2 + (b-a)\la_2)) + O((b-a) \la_2^{1/5} \la_3^{1/5}).  
\end{equation}
\end{lem} 
  
In order to get more accurate estimate for oscillatory sums, 
notice that if the values $\alpha$ и $\beta$ of $f'$ at the points $a$ и $b$ 
are located not too much close to integers, then the above estimate  
is more precise in first term. 
This estimate was established by Min in case of algebraic function~$f$. 
We will use the following more general result due to Liu (\cite{Liu}, Theorem~1). 
  
\begin{lem}\label{LemLiu} 
Let $f(y) \in C^5([a,b])$ be a real function, ${f''(y) > 0}$, 
and $R$, $U$ и $C_k$, ${k = 1,\ldots,6}$, be positive constants such that
\begin{gather}
	C_1R^{-1} \le |f''(y)| \le C_2R^{-1}, \qquad 
	|\beta_k(y)| \le C_kU^{2-k}, \quad U \ge 1, \quad 3 \le k \le 5, 
	\\
	\beta_k(y) = f^{(k)}(y)/f''(y). 
\end{gather}
Suppose $|3\beta_4(y) - 5\beta_3^2(y)| \ge C_6U^{-2}$ for any ${y \in [a,b]}$;
then
\begin{equation}
	\sum_{a \le y \le b} e^{2\pi i f(y)} = 
	\sum_{\alpha < k < \beta} |f''(y_k)|^{-1/2} e^{2\pi i\, (f(y_k) - ky_k + 1/8)} + \cE,  
\end{equation}
moreover,
\begin{equation}
\cE = \cE_1 + \cE_2 + 
	 O(\ln(2 + (b-a)R^{-1})) + O((b-a+R) U^{-1}) + 
	 O(\min\{ \sqrt{R}, \max\left(\frac1{\langle\alpha\rangle},\frac1{\langle\beta\rangle}\right) \}),  
\end{equation}
where 
\begin{gather}
	\cE_1 = b_\alpha |f''(y_\alpha)|^{-1/2} e^{2\pi i\,(f(x_\alpha) - \alpha x_\alpha + 1/8)}, 
	\\
	\cE_2 = b_\beta |f''(y_\beta)|^{-1/2} e^{2\pi i\,(f(x_\beta) - \beta x_\beta + 1/8)}, 
\end{gather}
$b_\alpha = 1/2$ if $\alpha \in \Set{Z}$, otherwise ${b_\alpha = 0}$, 
and $b_\beta$ is defined in the similar way; 
${\langle\alpha\rangle = \beta - \alpha}$ if ${\alpha \in \Set{Z}}$ and 
\begin{equation}
	\langle\alpha\rangle = \min_{n \in \Set{Z}} |n - \alpha|, \qquad
	\langle\beta\rangle = \min_{n \in \Set{Z}} |n - \beta| 
\end{equation}
if $\alpha \not\in \Set{Z}$ (respectively, $\beta \not\in \Set{Z}$). 
\end{lem}

\section[Free particle with exponential Hamiltonian]
{Free quantum particle with Hamiltonian $H(p) = E_0 + A e^{\eps p}$} 
  
\subsection{Outline of the method}\label{Outline} 
  
Consider a free quantum particle moving along the real line 
according to classical Hamiltonian 
${H(p) = \frac{p^2}{2}}$ (where $p$ is impulse variable). 
It generates an action $R_t \Maps \Set{R} \to \Set{R}$ of multiplicative group 
${\Set{R}_+ = \{t > 0\}}$ as follows. 
Solving equation 
\begin{equation}
	t\,\frac{\partial}{\partial p}H(p) = k, 
\end{equation}
we get
${p = k/t}$, and define ${R_t \Maps p \mapsto t^{-1} p}$.
Actually $R_t$ controls the dynamics of the stationary phase set 
for the sum 
\begin{equation}
	\cP(t) = \sum_{p \:\in\: [0,p^*] \,\cap\, \delta\Set{Z}} e^{itH(p)}, \qquad \delta > 0. 
\end{equation}
%
%
Obviously 
acting on $p^*$-periodic functions 
the map $R_{2}$ is dual to the map 
${\Set{T}_{p^*} \to \Set{T}_{p^*} \Maps p \to 2p \pmod{p^*}}$ 
which is a kind of {\it hyperbolic\/} transformation, 
where ${\Set{T}_{p^*} = \Set{R}/p^*\Set{Z}}$. This is the point which 
complicates the behavior of the trigonometric sum~$\cP(t)$.

The idea is to change the dynamics of the stationary phase set just by a small 
alteration of Hamiltonian adding terms for higher derivatives. 
This leads to exponential Hamiltonian 
\begin{equation}
	\Tilde H(p) = E_0 + \frac1{\eps^2} e^{\eps p} 
\end{equation}
generating action ${\Tilde R_t(p) = p - \ln t}$. 
Observe that (cf.\ \cite{Scardicchio}) 
\begin{equation}
	\Tilde H(p) = E_0 + \frac1{\eps^2} + \frac{p}{\eps} + \frac{p^2}{2} + \eps \frac{p^3}{6} + \ldots 
\end{equation}
Using Van der Corput's method (see \ref{sVanDerCorput}), we get
\begin{equation}
	S(t) \approx \int \Tilde R_t \phi \cdot \cL_t, 
\end{equation}
where 
${\Tilde R_t \phi(p) = \phi(\Tilde R_t(p))}$~, $\phi(p)$ is a test function  supported on $[a,b]$
and $\cL_t$~is generalized Legendre transform given by a discrete distribution on~$\Set{R}$.
Again using Van der Corput's method, we see that the sum $S(t)$ is flat.

\subsection{Main construction} 
\label{oscsummainconstr}
  
Consider trigonometric sum 
\begin{equation}
	S(t,\eps,q) = \frac1{\sqrt{q}} \sum_{y=0}^{q-1} e^{2\pi i\,t\omega(y)}, \qquad 
	\omega(y) = \omega_0 + \frac{q}{\eps^2} e^{\eps y/q}, 
\end{equation}
where $q \in \Set{N}$, ${0 < \eps < 1}$, ${\eps^{-1} \in \Set{N}}$. 
%
%
Adding a constant to $\omega(y)$ has no effect on $|\cP(t)|$, 
and during this section we assume that ${\omega_0 = 0}$. 
Our usual assumption on correlation of $t$ and $q$ will be ${1 \ll t \ll q^{1/3}}$. 
Let us list several simple properties of the function $\omega(y)$: 
\begin{itemize}
	\item[(a)] $\omega^{(r)}(y)$ increases for any ~$r$;
	\item[(b)] $\omega'(y) \ge \eps^{-1}$;
	\item[(c)] $\omega''(y) = \frac1{q}(1+O(\eps))$, 
		where large ``$O$'' is considered relatively to ${\eps \to 0}$, where ${0 \le y < q}$. 
\end{itemize}
Let $\cK(t)$ be the set of indexes $k$ such that $y_k(t)$ appears in $(0,q)$,
\begin{equation}
	\cK(t) = \left\{ k \in \Set{Z} \where y_k(t) \in (0,q) \right\}. 
\end{equation}
The set $\cK(t)$ is an integer segment, 
\begin{gather}
	\cK(t) = [k_0(t), k_1(t)], 
	\\ 
	K_0(t) = \frac{t}{\eps}, \quad 
	K_1(t) = \frac{t}{\eps}\, e^{\eps}, \quad
	\\ 
	k_0(t) = \min\{ k>K_0(t) \}, \quad 
	k_1(t) = \max\{ k<K_1(t) \}. 
\end{gather}

\begin{lem}\label{PAndPhi} 
Let $S(t,\eps, q)$ be the oscillatory sum defined 
in \ref {oscsummainconstr}, and assume that 
${\eps^{-1} \in \Set{N}}$, 
then $S(t,\eps, q)$ can be approximated by a 
sum over the set of stationary phases $\{y_k\}$ 
\begin{equation}
	S(t,\eps,q) = 
	\sum_{K_0(t) < k < K_1(t)} 
  	\frac{e^{2\pi i\, ( - ky_k)}}{\sqrt{t}} \gamma(y_k) + \frac{\cE(t)}{\sqrt{q}}, 
  \qquad
  \gamma(y) = e^{-\frac12 \eps y/q}, 
\end{equation}
where $\cE(t)$ is estimated as follows 
\begin{gather}
	\left|\frac{\cE(t)}{\sqrt{q}}\right| = 
	  O\left( \frac{\max\{t^{-1}, \ln t\}}{\sqrt{q}} \right) + 
	  O\left(
	   \frac1{\sqrt{q}} 
	   \min\{ 
	    \sqrt{\frac{q}{t}}, \nu_\eps(t)
	   \}
	   \right),  
	\\
	\nu_{\eps}(t) = 
	    \frac1{ \left\|\frac{t}{\eps}\right\|_{\Set{Z}} } + 
	    \frac1{ \left\|\frac{t}{\eps}e^{\eps}\right\|_{\Set{Z}} }, 
\end{gather}
and $\|x\|_{\Set{Z}} = \min_{n \in \Set{Z}} |n - x|$. 
For any $t$ the sum $S(t,\eps,q)$ satisfies the evident condition ${|S(t,\eps,q)| \le \sqrt{q}}$. 
\end{lem}

\begin{rem}
Investigating the sum $S(t,\eps,q)$ small and large O's are considered with respect 
to ${\eps \to 0}$, ${t \to \infty}$ and ${q \to \infty}$. 
\end{rem}

\begin{rem}
We will apply this lemma to find upper bounds for the error term 
with respect to different integral norms. So we can skip 
some countable set of points $t$ in our considerations. 
\end{rem}

\begin{proof} 
Let $f(y) = t\omega(y)$ and note that ${\omega'(y) = \eps q^{-1} \omega(y)}$ 
and ${f'(y) = \eps q^{-1} f(y)}$. Applying the main Lemma \ref{LemLiu} 
with ${a = 0}$, ${b = q}$, we have
\begin{multline}
	S(t,\eps,q) = 
		\sum_{K_0(t) < k < K_1(t)} 
  		\frac{e^{2\pi i\, (f(y_k) - ky_k)}}{\sqrt{q\,|f''(y_k)|}} + q^{-1/2}\cE(t) = \\
	=
	\sum_{K_0(t) < k < K_1(t)} 
  	\frac{e^{2\pi i\, (f(y_k) - ky_k)}}{\sqrt{t}} \gamma(y_k) + q^{-1/2}\cE(t) 
  = 
	\sum_{K_0(t) < k < K_1(t)} 
  	\frac{e^{2\pi i\, (- ky_k)}}{\sqrt{t}} \gamma(y_k) + q^{-1/2}\cE(t), 
\end{multline}
since  $f(y_k) = \eps^{-1} q f'(y_k) = \eps^{-1} q k \in \Set{Z}$,  
\begin{equation}
	q\,|f''(y_k)| = q\, \frac{t}{q} e^{\eps y/q} = t \gamma^{-2}(y_k), \qquad 
	\gamma(y) = e^{-\frac12 \eps y/q}, 
	\qquad
	e^{\eps y/q} = 1 + O(\eps), 
\end{equation}
and $\cE(t)$ has the following form (notation from Lemma \ref{LemLiu}) 
\begin{multline}
	\cE(t) = \cE_1(t) + \cE_2(t) + \\ +
	 O(\ln(2 + (b-a)R^{-1})) + O((b-a+R) U^{-1}) + 
	 O(\min\{ \sqrt{R}, \max\left(\frac1{\langle\alpha\rangle},\frac1{\langle\beta\rangle}\right) \}),  
\end{multline}
where $a=0$, $b=q$ and values $R$ and $U$ are calculated below: 
\begin{gather}
	b-a = q, \qquad K_1(t) - K_0(t) = t+O(t\eps) \sim t, 
	\\
	\frac{t}{q} \le f''(y) \le \frac{t}{q}(1+O(\eps)), \qquad R = \frac{q}{t}, 
	\\
	f^{(2+j)}(y) = \frac{\eps^j}{q^j} f''(y), \qquad 
	\beta_{2+j} = \frac{\eps^j}{q^j}, \qquad j = 1,2,3, 
	\\
	U = \frac{q}{\eps} > 1, 
	\\
	|3\beta_4 - 5\beta_3^2| = \left|3\frac{\eps^2}{q^2} - 5\left(\frac{\eps}{q}\right)^2\right| 
		= 2\frac{\eps^2}{q^2} = 2U^{-2}. 
\end{gather}
So, applying Lemma \ref{LemLiu} we obtain
\begin{multline}
	\cE(t) = 
	\cE_1(t) + \cE_2(t) + 
	 O(\ln(2 + t)) + O\left((q+t^{-1}q) \frac{\eps}{q}\right) + \cE_3(t) = \\ 
	=
	\cE_1(t) + \cE_2(t) + 
	 O(\max\{t^{-1},\ln(t)\}) +  \cE_3(t), 
\end{multline} 
where
\begin{equation}
	 \cE_3(t) = 
	  O(
	   \min\{ 
	    \sqrt{\frac{q}{t}}, 
	    \max\left(\frac1{\langle\alpha\rangle},\frac1{\langle\beta\rangle}\right) 
	   \}
	   ). 
\end{equation}
The terms $\cE_1(t)$ and $\cE_2(t)$ are supported on a countable set. 
To compute $\cE_3(t)$ look at the set 
\begin{equation}
	\Lambda(t)=\{y_k(t)\}, \quad y_k(t)= \frac{q}{\eps}\ln \frac{\eps k}{t}. 
\end{equation}
The set $\Lambda(t)$ is moving rigidly along the line when $t$ changes .  
The term ${\cE_3(t) = O(\sqrt{q/t})}$ is large whenever some point in $\Lambda(t)$ 
is close to border of $[0,q]$.
At the same time for a typical $t$ we have ${\cE_3(t) = O(1)}$.
Further, we have the following explicit representation for $\langle\alpha\rangle$ and $\langle\beta\rangle$: 
\begin{figure}[th]
	\begin{center}
		\includegraphics[width=100mm]{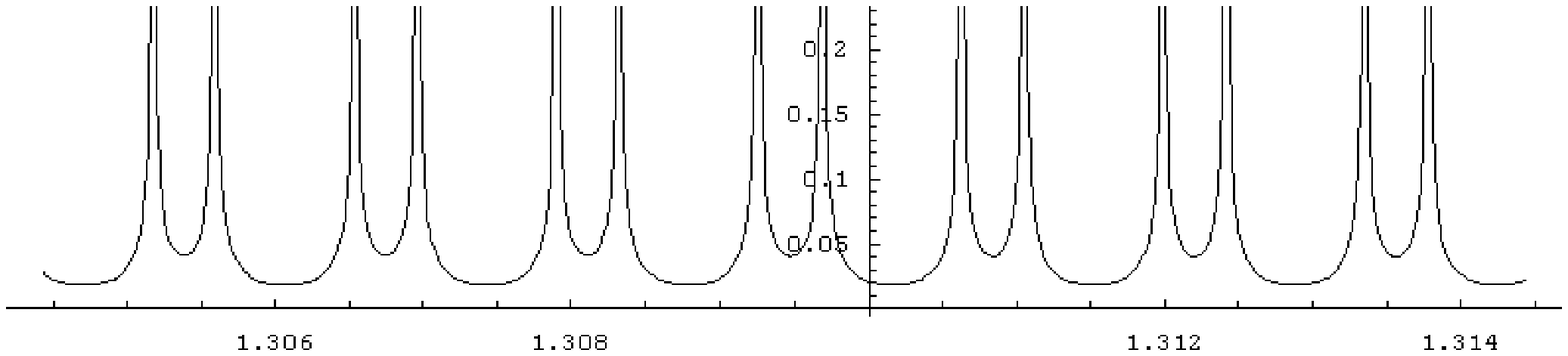}
	\end{center}
	\caption{Estimating $\cE_3(t)$} 
\end{figure}
\begin{equation}
	\langle\alpha\rangle = O( \|f'(0)\|_{\Set{Z}} ) = O\left( \left\|\frac{t}{\eps}\right\|_{\Set{Z}} \right), 
	\qquad 
	\|x\|_{\Set{Z}} = \min_{n \in \Set{Z}} |n - x|. 
\end{equation}
Similarly,
\begin{equation}
	\langle\beta\rangle = O\left( \left\|\frac{t}{\eps}e^{\eps}\right\|_{\Set{Z}} \right). 
\end{equation}
Thus, omitting the terms $\cE_1(t)$ and $\cE_2(t)$ with countable supports we have 
\begin{equation}
	\left|\frac{\cE(t)}{\sqrt{q}}\right| = 
	  O\left( \frac{\max\{t^{-1}, \ln t\}}{\sqrt{q}} \right) + 
	  O\left(
	   \frac1{\sqrt{q}} 
	   \min\{ 
	    \sqrt{\frac{q}{t}}, \nu_\eps(t)
	   \}
	   \right). 
\end{equation}
\end{proof}

\begin{lem}\label{PAndPhiIntegral} 
Given $\eps$ and $q$ consider the sum $S(t,\eps,q)$ for ${t \in [t_1,t_2]}$, 
and suppose that all conditions of the previous lemma are satisfied 
as well as the following requirements 
\begin{equation}
	t_1 \ge 1, \quad t_2-t_1 \ge \eps, \quad t_2 \le q, \quad \text{and} \quad \eps^{-1} \le q. 
\end{equation}
Then 
\begin{gather}
	\|q^{-1/2}\cE|_{[t_1,t_2]}\|_\infty = 
		O\left( \frac{\max\{1, \ln t_2\}}{\sqrt{q}} \right) + 
	  O\left( \frac1{\sqrt{t_1}} \right), 
	\\ 
	\|q^{-1/2}\cE|_{[t_1,t_2]}\|_1 \lesssim (t_2-t_1)\, \frac{\ln q}{\sqrt{q}}. 
\end{gather}
\end{lem}

\begin{proof}
Since ${t \ge \eps}$
\begin{equation}
	\int_{t_1}^{t_2} |\cE(t)| \,dt \lesssim 
		\frac{t_2-t_1}{\eps} \cdot \int_0^{\eps} \min\{\sqrt{\frac{q}{t}}, \frac{\eps}{\tau}\} \,d\tau 
	\lesssim (t_2-t_1)\ln q, 
\end{equation}
hence,
\begin{equation}
	\int_{t_1}^{t_2} \left|\frac{\cE(t)}{\sqrt{q}}\right| \,dt \lesssim 
	\frac{(t_2-t_1)\ln t_2 + (t_2-t_1)\ln q}{\sqrt{q}} \lesssim 
	(t_2-t_1)\, \frac{\ln q}{\sqrt{q}}. 
\end{equation}
\end{proof}

Our purpose is to apply this series of lemmas around van der Corputs's method 
to prove that ${|S(t,\eps,q)| \approx 1}$ for most points $t$. Therefore 
let us formulate the following improved lemma.

\begin{lem}
In the scope of Lemma \ref{PAndPhiIntegral} the expression for $S(t,\eps,q)$ 
can be simplified as follows:  
\begin{equation}
	S(t,\eps,q) = \frac{1}{\sqrt{t}} 
	\sum_{K_0(t) < k < K_1(t)} 
	  e^{2\pi i\, ( - ky_k)} + O(\eps \sqrt{t}) + \frac{\cE(t)}{\sqrt{q}}, 
\end{equation}
where $\frac{\cE(t)}{\sqrt{q}}$ is estimated as shown in Lemmas \ref{PAndPhi} and \ref{PAndPhiIntegral}. 
\end{lem}

\begin{proof}
The lemma follows from the statement of lemma \ref{PAndPhiIntegral} and 
simple observation that the sum over $k$ contains no more than 
${K_1(t)-K_0(t) \sim t}$ terms and the fact that $\gamma(y)$ is $\eps$-close to $1$ 
for ${y \in [0,q]}$, namely, ${\gamma(y) = e^{-\frac12 \eps y/q} = 1 + O(\eps)}$. 
\end{proof}

At this point noitce that in the expression for $S(t,\eps,q)$ we see 
a reduced (with smaller number of terms) oscillatory sum of the same form as $S(t,\eps,q)$, 
and we can iterate the reduction procedure applying van der Corput's method to that sum. 
Denoting the new frequency function as $\eta_{\eps,q}(k)$,
\begin{equation}
	\eta_{\eps,q}(k) = -ky_k, 
\end{equation}
let us see how it looks like for the sequence of stationary phases $y_k$ 
generated by an exponential Hamiltonian $\omega(y)$: 
\begin{equation}
	\eta_{\eps,q}(k) = -ky_k = 
		-\frac{qk}{\eps}\ln\frac{\eps k}{t} =
 		\frac{q(\ln t-\ln\eps)}{\eps}\cdot k - q\cdot \frac{k \ln k}{\eps} = 
 		x_{\eps,q}(t)k - q \cdot \Omega_{\eps}(k), 
\end{equation} 
where 
\begin{equation}
	x_{\eps,q}(t) = \frac{q(\ln t-\ln\eps)}{\eps} 
\end{equation}
and
\begin{equation}
	\Omega_{\eps}(k) = \frac1{\eps} k \ln k. 
\end{equation}

It appears that starting from the property of constant distance between 
adjacent stationary phases we can find a variety of functions $\omega(y)$ 
(having continuum cardinality) generating flat reduced sums. 
At the same time we will see that surprisingly there are infinite subsequence 
of flat sums in the sequence of frequence functions $\eta_{\eps,q}(k)$ 
without any modification in the construction. To see that we observe 
that modulo one ${Y(q) = q \cdot \Omega_{\eps}(k)}$ is a point 
in the torus $\prod_k \Set{T}$ ruled by the dynamics of rigid shift 
by the vector $\Omega_{\eps}(k)$.

\begin{lem}
Let $k$ range over a finite set $\cK$. 
For any ${\delta > 0}$ there exist infinitely many ${q = q_j \to \infty}$ 
such that for each ${k \le m\eps^{-1}}$ 
\begin{equation}
	q_j \cdot \Omega_{\eps}(k) = -\Omega_{\eps}(k) + O(\delta) \pmod{1}. 
\end{equation}
\end{lem} 
  
\begin{proof}
The statement follows from Poincar\'e recurrence theorem 
for torus shift on $\prod_{k \in \cK}\Set{T}$ by the vector $\Omega_{\eps}(k)$. 
\end{proof} 

\begin{cor}
It is important that if $t$ ranges over a bounded segment, i.e.\ ${t \in [t_1,t_2]}$, 
then $k$ ranges over a finite set ${[K_0(t_1),K_1(t_2)] \cap \Set{Z}}$. 
Thus, the statement of the previous lemma concerning the dynamics of $Y(q)$ 
is true for a sequence $q_j$ simultaneously for all integer $k$ that 
can appear in $[K_0(t),K_1(t)]$ for all ${t \in [t_1,t_2]}$. 
\end{cor}

\begin{lem} 
The oscillatory sum defined by $\Omega_{\eps}(k)$ 
approximates the reduced sum given by the function $\eta_{\eps,q}(k)$. If 
\begin{equation}
	q \cdot \Omega_{\eps}(k) = -\Omega_{\eps}(k) + O(\delta) \pmod{1} 
\end{equation}
then 
\begin{equation}
	\left|
		\frac1{\sqrt{t}} \sum_{K_0(t) < k < K_1(t)}
			e^{2\pi i\,\eta_{\eps,q}(k)}
		-
		\frac1{\sqrt{t}} \sum_{K_0(t) < k < K_1(t)}
			e^{2\pi i\,(x_{\eps,q}(t)k - \Omega_{\eps}(k))}
	\right|
	= O(\delta\sqrt{t}). 
\end{equation}
\end{lem}

\begin{proof}
We know that $\eta_{\eps,q}(k) = x_{\eps,q}(t)k - q \, \Omega_{\eps}(k)$ and 
the lemma follows from the estimate of the total number of terms in the sum, 
${K_1(t)-K_0(t) \sim t}$. 
\end{proof}

Now let us study the oscillatory sum defined by $\Omega_{\eps}(k)$ 
with additional linear term $xk$ 
\begin{equation}
	\Phi(t,x,\eps) = \frac1{\sqrt{t}} 
	 \sum_{K_0(t) < k < K_1(t)} \exp(2\pi i\,(xk + \Omega_{\eps}(k))). 
\end{equation}

\begin{rem}
The fundamental observation concerning $\Omega_{\eps}(k)$ consists of the following. 
When $k$ passes the interval $(K_0(t),K_1(t))$, 
\begin{equation}
	k \in (K_0(t),K_1(t)) = \left( \frac{t}{\eps}, \frac{t}{\eps}e^\eps \right) 
\end{equation}
the derivative ${\Omega_{\eps}(k) = \frac1{\eps}(1+\ln k)}$ increases exactly by one:
\begin{equation}
	\Omega_\eps'(K_1(t)) - \Omega_\eps'(K_0(t)) = 
	\frac1{\eps}\ln \frac{t}{\eps} - \frac1{\eps}\ln\left( \frac{t}{\eps}e^{\eps} \right) = 1, 
\end{equation}
hence, we get (almost surely) exactly one stationary phase ${k^* \in (K_0(t),K_1(t))}$. 
Also let us point that the terms of the sum do not depend on $t$, just $K_0$ and $K_1$ do.
\end{rem}

\begin{lem}\label{LemPhi}
Let the following conditions be satisfied: 
\begin{gather}
	t \ge 1, \qquad 
	0 < \eps < 1, \qquad \eps^{-1} \in \Set{N}, 
\end{gather}
and let $(k^*,\ell)$ be a unique solution of equation 
${x+\Omega_{\eps}'(k^*) = \ell}$, ${\ell \in \Set{Z}}$, ${k^* \in (K_0(t),K_1(t))}$ (if exists). Then 
\begin{gather}
	\Phi(t,x,\eps) = e^{2\pi i\, (\Omega_\eps(k^*) + x k^* - \ell k^*)} + t^{-1/2}\cE_\Phi(t,x), 
	\\
	t^{-1/2}\cE_\Phi(t,x) = O(\eps) + 
	 \frac1{\sqrt{t}} \,O\!\left( \min\{ \sqrt{t}, \frac1{\| \frac1{\eps}\ln\frac{t}{\eps} + x \|_{\Set{Z}}} \} \right). 
\end{gather}
\end{lem}

\begin{proof} 
Let us apply again Lemma \ref{LemLiu} to ${xk + \Omega_{\eps}(k)}$. Calculating 
\begin{equation}
	(xk + \Omega_{\eps}(k))'_k = x + \frac1{\eps}(1 + \ln k), \qquad 
	(xk + \Omega_{\eps}(k))''_k = \frac1{\eps k}, 
\end{equation}
notice that $k \sim t/\eps$ and we can take 
\begin{equation}
	R^{-1} = \frac1{\eps \, \frac{t}{\eps}}, \quad \text{and} \quad R = t. 
\end{equation} 
Further, let us define 
\begin{gather}
	\beta_3(k) = -\frac1{k}, \quad 
	\beta_4(k) = \frac2{k^2}, \quad 
	\beta_5(k) = -\frac6{k^3}, 
	\\
	U = \frac{t}{\eps} \sim k, \qquad 
	|3\beta_4 - 5\beta_3^2| = \frac1{k^2} = U^{-2}, 
\end{gather}
Recall that ${\Omega_\eps'(K_1(t)) - \Omega_\eps'(K_0(t)) = 1}$. So, 
we the following representation for the sum $\Phi(t,x,\eps)$:
\begin{equation}
	\Phi(t,x,\eps) = 
		\frac1{\sqrt{t\,|\Omega_{\eps}''(k^*)|}} 
		e^{2\pi i\, (\Omega_\eps(k^*) + x k^* - \ell k^*)} + q^{-1/2} \cE_0(t,x), 
\end{equation}
where $(k^*,\ell)$ are found as a unique solution (if exists) of the equation 
${x+\Omega_{\eps}'(k^*) = \ell}$ with ${\ell \in \Set{Z}}$, ${k^* \in (K_0(t),K_1(t))}$. 
Clearly the coefficient near exponent in the right-hand side is close to $1$ in absolute value, 
\begin{equation}
	\frac1{\sqrt{t\,|\Omega_{\eps}''(k^*)|}} = 
	\left( t\,\frac1{\eps k^*} \right)^{-1/2} = 
	\left( t\,\frac1{\eps \cdot \eps^{-1}t(1+O(\eps))} \right)^{-1/2} = 1 + O(\eps), 
\end{equation}
and using Lemma \ref{LemLiu} we can compute the following estimate for the error term $\cE_0$ 
\begin{equation}
	\cE_0(t,x) = \cE_1 + \cE_2 + O(1) + \cE_3, 
\qquad
	\cE_3 = O( \min\{ \sqrt{t}, \tilde\nu \} ), 
\end{equation}
with 
\begin{equation}
	\tilde\nu = 
	\frac1{\| K_0(t) + x \|_{\Set{Z}}} = 
	\frac1{\| \frac1{\eps}\ln\frac{t}{\eps} + x \|_{\Set{Z}}}. 
\end{equation}
Finally, we have the following expression for the common error term 
\begin{equation}
	t^{-1/2}\cE_\Phi(t,x) = O(\eps) + 
		\frac1{\sqrt{t}} \,
		O\!\left( 
			\min\{ \sqrt{t}, \frac1{\| \frac1{\eps}\ln\frac{t}{\eps} + x \|_{\Set{Z}}} \} 
		\right). 
\end{equation}
%
\end{proof}

Now to join the reduction procedure from $S(t,\eps,q)$ to $\Phi(t,x,\eps)$ and 
approximation given by Lemma \ref{LemPhi} we have to calculate boundary effect term 
${\frac1{\eps}\ln\frac{t}{\eps} + x_{\eps,q}(t)}$ as well as the argument 
of the exponent, namely, ${\Omega_\eps(k^*) + x_{\eps,q}(t) k^* - \ell k^*}$ 
(which is of independent interest but does not influence to flatness). 
So,
\begin{equation}
	v_{\eps,q}(t) = \frac1{\eps}\ln\frac{t}{\eps} + x_{\eps,q}(t) = 
	(q-1)\frac1{\eps}\ln\frac{t}{\eps}. 
\end{equation}
Further, remark that 
\begin{equation}
	\Omega_\eps(k^*) + x_{\eps,q}(t) k^* - \ell k^* = k^*(\ell - \frac1{\eps}) 
\end{equation}
do not depend ``directly'' on $x_{\eps,q}$ but via $k^*(t)$ and $\ell(t)$. 
The rate of grows of this doubly reduced frequency function can be estimated as follows 
\begin{equation}
	\Omega_\eps(k^*) + x_{\eps,q}(t) k^* - \ell k^* = O\bigl(k^*(t)\ell(t)\bigr) 
	\sim \frac1{\eps} k^*(t) \ln k^*(t) 
	\sim \frac{t}{\eps^2} \ln \frac{t}{\eps}. 
\end{equation}

\begin{lem}\label{PAndPhiComplex} 
In the scope of Lemmas \ref{PAndPhi} and \ref{LemPhi} 
for an arbitrary ${\delta > 0}$ there exist a sequence 
${q = q_j \to \infty}$ such that for ${t \ge 1}$ 
\begin{equation}
	S(t,\eps,q) = \exp\left(2\pi i\, k^*(t)\bigl(\ell(t) - \eps^{-1}\bigr)\right) + \tilde\cE(t)
\end{equation}
with an error term (defined up to countable set of points $t$) 
\begin{multline}
	\tilde\cE(t) = 
		O(\delta \sqrt{t}) + 
		O(\eps \sqrt{t}) + 
		O\left( \frac{\max\{1, \ln t\}}{\sqrt{q}} \right) + 
	  O\left(
	   \frac1{\sqrt{q}} 
	   \min\{ 
	    \sqrt{\frac{q}{t}}, \nu_\eps(t)
	   \}
	   \right) + \\
	  +
		\frac1{\sqrt{t}} \,
		O\!\left( 
			\min\{ \sqrt{t}, \frac1{\| (q-1)\frac1{\eps}\ln\frac{t}{\eps} \|_{\Set{Z}}} \} 
		\right). 
\end{multline}
\end{lem}

\begin{proof} 
The proof is just a compilation of Lemmas \ref{PAndPhi} and \ref{LemPhi}. 
\end{proof}

\begin{lem}\label{PAndPhiComplexInt} 
Preserving conditions of Lemma \ref{PAndPhiComplex} suppose that ${t \in [t_1,t_2]}$ and 
${t_2-t_1 \ge (q-1)^{-1}\eps t_2}$. Then 
\begin{equation}
	\Bigl\| 
		\bigl| S(t,\eps,q)|_{[t_1,t_2]} \bigr| - 1 
	\Bigl\|_1
	= 
	O\left( (t_2-t_1)\, \sqrt{t_2}\, \max\{\delta,\eps\} \right) 
	+
	O\left( (t_2-t_1)\, \frac{\ln q}{\sqrt{q}} \right) 
	+
	O\left( (t_2-t_1)\, \frac{\ln t_2}{\sqrt{t_1}} \right) 
	.  
\end{equation}
\end{lem}

\begin{proof} 
We have to estimate the influence of the last term in the error given by the previous lemma, 
\begin{equation}
	\int_{t_1}^{t_2} 
		\frac1{\sqrt{t}} 
		\min\{ \sqrt{t}, \frac1{\| (q-1)\frac1{\eps}\ln\frac{t}{\eps} \|_{\Set{Z}}} \} \,dt 
	\lesssim (t_2 - t_1)\, \frac{\ln t_2}{t_1}. 
\end{equation}
\end{proof}

\begin{lem}
If the parameters $t_1,t_2,\eps,q$ satisfy the following requirements 
\begin{equation}
	t_1 \ge 1, \quad t_2-t_1 \ge \eps, \quad t_2 \le \frac1{\eps} \le \frac{q}{\ln q}, \qquad 
	\delta \le \eps, 
\end{equation}
then for a sequence ${q_j \to \infty}$ 
\begin{gather}
	\Bigl\| 
		\bigl| S(t,\eps,q_j)|_{[t_1,t_2]} \bigr|^2 - 1 
	\Bigl\|_1 
	\lesssim 
	\Bigl\| 
		\bigl| S(t,\eps,q_j)|_{[t_1,t_2]} \bigr| - 1 
	\Bigl\|_1
	= 
	O\left( (t_2-t_1)\, \frac{\ln t_2}{\sqrt{t_1}} \right), 
\end{gather}
\end{lem}

\begin{proof}
To established the improved estimate of the error term we just have to 
observe that all the components of the estimate are dominated by the term 
$O\left( (t_2-t_1)\, t_1^{-1/2}\ln t_2 \right)$. 
The first inequality of the lemma follows from the fact that ${\tilde\cE(t) = O(1)}$. 
\end{proof}

Recall that the initial frequency function $\omega(y)$ satisfies 
differential equation ${\omega'' = a\omega'}$ with solution 
${\omega(y) = E_0 + e^{a(y-y_0)}}$ for arbitrary $E_0,y_0$. 
By this point we studied $\omega$ not taking into account 
correlation of additional multiplier ${m = e^{-ay_0}}$ near $t$ 
and the length $q$ of the segment $[0,q]$.

\begin{lem}\label{LemSmInt}
Let us consider the sum 
\begin{equation}
	S_m(\tau,\eps,q) = 
		\frac1{\sqrt{q}} \sum_{y=0}^{q-1} e^{2\pi i\, \tau\omega(y)} 
\end{equation}
with 
\begin{equation}
	\omega(y,m) = m\frac{q}{\eps^2} e^{\eps y/q}. 
\end{equation}
Suppose that ${t_1 = m\tau_1}$, ${t_2 = m\tau_2}$ as well as $\eps$ and $q$ 
meet the requirements of the previous lemma. 
Then for a sequence ${q_j \to \infty}$ 
\begin{gather}
	\Bigl\| 
		\bigl| S(t,\eps,q_j)|_{[\tau_1,\tau_2]} \bigr|^2 - 1 
	\Bigl\|_1 
	= 
	O\left( (\tau_2-\tau_1)\, \frac{\ln(m\tau_2)}{\sqrt{m\tau_1}} \right). 
\end{gather}
\end{lem}

Applying this lemma in the situation when $\tau_1$ and $\tau_2$ are fixed 
and $m$ is sufficiently large we come to the following theorem.

\begin{thm}
For given $0 < a < b$, ${\eps > 0}$ and ${\delta > 0}$ 
there exists $m_0$ such that for any ${m \ge m_0}$ 
there exists an infinite sequence $q_j$ generating 
trigonometric sums with exponential staircase frequency function
\begin{equation}
	\omega(y,m) = m\frac{q_j}{\eps^2}e^{\eps y/q_j}
\end{equation}
which are $\delta$-flat in $L^1([a,b])$.
\end{thm}

%% file: proof_riesz.tex
\section[Convergence of Riesz products]{Convergence of Riesz~products for exponential staircase flow}

\subsection{Preliminary lemmas}

Let $R(t)$ be the correlation function of a measurable flow $T^t$. 
Observe that if the flow atcs on a Lebesgue space then it obeys 
the following continuity property. If ${f,g \in L^2(X)}$ then 
\begin{equation}
	\scpr<T^f,g> \to \scpr<f,g> \quad \text{as} \quad t \to 0. 
\end{equation}

\begin{lem}
Suppose that $R(t)$ is characteristic function of a spectral measure $\sigma_f$, i.e.\ 
\begin{equation}
	R(t) = \hat\sigma_f(t) = \int_{\Set{R}} e^{2\pi i\, t x} \,d\sigma_f(x). 
\end{equation}
Let $\nu_n$ be a sequence of finite positive measures and let ${R_n = \hat\nu_n}$. 
Then $\nu_n$ converges weakly to $\sigma_f$ iff ${R_n(t) \to R(t)}$ for any~$t$. 
\end{lem}

\begin{proof}
The statement follows from continuity of $R(t)$ at point ${t = 0}$. 
\end{proof}

Throughout this section let $\nu_n$ be a sequence of finite positive measures. 
Evidently if we have to prove weak convergence of $\nu_n$ 
it is sufficient to watch $\int_{\Set{R}} \psi \,d\nu_n$ for a dense set 
of functions $\psi$ in $C_0(\Set{R})$, 
for example a dense set of smooth functions with bounded support. 
Let us denote $C_0(a,b)$ the set of continuous functions $\psi(t)$ on $(a,b)$ 
such that ${\psi(t) \to 0}$ if ${t \to a}$ or ${t \to b}$. 

\begin{lem}
Let $\Phi(t)$ be a function such that it is locally $L^1$ and 
for any function ${\psi \in C_0(a,b)}$ with ${0 < a < b}$ (or ${a < b < 0}$) 
\begin{equation}
	\int \psi \,d\nu_n \to \int \psi(t) \,\Phi(t) \,dt. 
\end{equation}
Suppose that ${\nu_n \to \eta}$ weakly. 
Then $\Psi \in L^1(\Set{R}$ and ${\|\Psi\|_1 \le \|\eta\|}$. 
\end{lem}

\begin{proof}
Let us denote 
\begin{equation}
	\scpr<\eta,\psi> = \int \psi \,d\eta. 
\end{equation}
By definition of weak convergence ${\scpr<\nu_n,\psi> \to \scpr<\eta,\psi>}$, and, so 
\begin{equation}
	\liminf_{n \to \infty} |\scpr<\nu_n,\psi>| \le |\scpr<\eta,\psi>| \le \|\eta\|\cdot \|\psi\|_\infty. 
\end{equation}
At the same time we know that ${\scpr<\nu_n,\psi> \to \scpr<\Phi,\psi>}$, hence, 
\begin{equation}
	|\scpr<\Phi,\psi>| \le \|\eta\|\cdot \|\psi\|_\infty, 
\end{equation}
and the statement follows from the equality ${\|\Phi\|_1 = \|\la_\Phi\|}$ 
where $\la_\Phi$ is the measure with density $\Phi$. 
\end{proof}

\begin{lem}
Let $\eta$ be a finite positive measure on the line. 
If ${\int \psi \,d\eta = \int \psi(t) \,\Phi(t) \,dt}$ for any ${\psi \in C_0(a,b)}$ 
then ${\eta = \eta(\{0\})\df_0 + \la_\Phi}$ where $\df_0$ is the unit atomar measure in $0$ 
and $\la_\Phi$ is the measure with density $\Phi$. 
\end{lem}

\begin{proof}
It is possible to check coincidence on positive continuous functions $\psi$ 
with bounded support since its linear combinations are dense in $C_0(\Set{R})$. 
For any ${\eps > 0}$ there exists $\alpha$ such that 
\begin{equation}
	\int_{-\alpha}^\alpha \Phi(t) \,dt < \eps 
\end{equation}
and 
\begin{equation}
	\eta|_{[-\alpha,\alpha]} = \eta(\{0\})\df_0 + \eta_{(\alpha)} 
	\quad \text{with} \quad 
	\|\eta_{(\alpha)}\| < \eps. 
\end{equation}
Let us split $\psi$ in two parts, ${\psi = \psi_0 + \psi_1}$ with the following properties: 
\begin{equation}
	\supp \psi_0 \subseteq [-\alpha,\alpha], \qquad 
	\supp \psi_1 \subseteq [-L,-\alpha/2] \cup [\alpha/2,L] 
\end{equation}
and ${0 \le \psi_0(t),\psi_1(t) \le \psi(t)}$ 
for any~$t$. 
We have 
\begin{equation}
	\scpr<\eta,\psi> = \scpr<\eta,\psi_0> + \scpr<\eta,\psi_1> = 
	\eta(\{0\})\psi(0) + \scpr<\eta_{(\alpha)},\psi_0> + 
	\scpr<\Phi,\psi_1-\psi> + \scpr<\Phi,\psi>. 
\end{equation}
Using inequality $|\psi(t)| \le \|\psi\|_\infty$, we can estimate small terms in the above formula: 
\begin{equation}
	|\scpr<\eta_{(\alpha)},\psi_0>| \le \eps \cdot \|\psi\|_\infty, \qquad 
	|\scpr<\Phi,\psi_1-\psi>| = |\scpr<\Phi,\psi_0>| < \eps \cdot \|\psi\|_\infty. 
\end{equation}
Since $\eps$ is arbitrary we get $\scpr<\eta,\psi> = \eta(\{0\})\psi(0) + \scpr<\Phi,\psi>$. 
\end{proof}

\begin{lem}\label{LemConvOfTrigSums}
Consider positive functions $f(t)$ and $Q_n(t)$ on the real line $\Set{R}$ 
satisfying the following conditions:
\begin{itemize}
	\item[(i)] $f \in L^1(\Set{R})$; 
	\item[(ii)] $Q_n$ are continuous;
	\item[(iii)] $\dstyle \int_{a_n}^{b_n} |Q_n(t)-1| \,dt < \eps_n^2$ and 
		$\dstyle \int_{-b_n}^{-a_n} |Q_n(t)-1| \,dt < \eps_n^2$, 
		where ${0 < a_n < b_n}$;
	\item[(iv)] intervals $(a_n,b_n)$ monotonously increase and exhaust $(0,+\infty)$, exactly, 
		${a_{n+1} < a_n}$, ${b_n < b_{n+1}}$, ${\lim_n a_n = 0}$ and ${\lim_n b_n = +\infty}$; 
	\item[(v)] $\dstyle \sum_{n=1}^\infty \eps_n < \infty$. 
\end{itemize} 
The infinite product of $f$ and $Q_n$ converges in $L^1$ on any interval $(a^*,b^*)$ or $(-b^*,-a^*)$, 
where ${0 < a^* < b^*}$. 
\end{lem}

\begin{proof}
Without lost of generality, passing from $f(t)$ to $\tilde f(t) = f Q_1 \cdots Q_{n_1}(t)$ 
we can assume that ${(a^*,b^*) \subseteq (a_1,b_1)}$. Let us consider sets 
\begin{equation}
	E_n = \bigl\{ t \in (a^*,b^*) \where |Q_n(t)-1| > \eps_n \bigr\}. 
\end{equation}
By Chebyshev inequality $\la (E_n) < \eps_n^2/\eps_n = \eps_n$. 
Using Borel--Cantelli lemma we can deduce from convergence ${\sum_n\la(E_n) < \infty}$ 
that almost surely a point ${t \in (a^*,b^*)}$ belongs to finitely many sets $E_n$, 
hence, the product ${f(t)\prod_n Q_n}$ converges by Lebesgue dominanted convergence theorem. 
Indeed, if we restrict functions to the set 
\begin{equation}
	A_N = (a^*,b^*) \sms \bigcup_{n > N} E_n, \qquad \la(A_n) \to b^*-a^*, 
\end{equation}
then by pointwise convergence 
\begin{equation}
	\prod_{n = N+1}^m Q_n(t) \to \prod_{n > N} Q_n(t), \qquad |Q_n(t)-1| < \eps_n \quad \text{if} \quad n > N, 
\end{equation}
we get convergence in $L^1$ 
\begin{equation}
	\Pi_m = f\prod_{n=1}^m Q_n \to f\prod_n Q_n 
\end{equation}
for the restrictions to interval $(a^*,b^*)$. 
\end{proof}

\begin{lem}\label{LemConvToLOne}
Preserving conditions of the previous lemma let us assume that 
the partial products $\Pi_m(t)$ are $L^1$-densities of measures on $\Set{R}$ 
which converge weakly to a finite positive measure~$\eta$, i.e.\ 
for any ${\psi \in C_0(\Set{R})}$ 
\begin{equation}
	\int \psi(t) \,\Pi_m(t) \,dt \to \scpr<\eta,\psi>. 
\end{equation}
Then there is a function $\Phi \in L^1(\Set{R})$ such that 
${\eta = \eta(\{0\})\df_0 + \la_\Phi}$ and for Lebesgue almost every ${t \in \Set{R}}$ 
we have ${\Phi(t) = f(t)\prod_n Q_n(t)}$. 
\end{lem}

\begin{proof}
Consider an interval $(a^*,b^*)$ with ${0 < a^* < b^*}$. 
By $L^1$-convergence $\Pi_n|_{(a^*,b^*)}$ converges weakly to the measure 
with density $f(t)\prod_n Q_n(t)$ restricted to $(a^*,b^*)$. 
This implies that the weak limit $\eta$ of $\Pi_n$ acting on functions 
${\psi \in C_0(a^*,b^*)}$ coincide with ${\Phi(t) = f(t)\prod_n Q_n(t)}$. 
Using above lemmas we found that ${\Phi \in L^1(\Set{R})}$ and 
${\eta = \eta(\{0\})\df_0 + \la_\Phi}$. 
\end{proof}

%


\subsection{Spectral type of exponential staircase flows}

Consider exponetial staircase construction of rank one flow which is characterized by 
the following parameters: numbers of subcolumns $q_n$ and staircase grade values $m_n$. 
We define $\eps_n$ by the equality 
\begin{equation}
	h_n = \frac{m_n}{\eps_n}, 
\end{equation}
and require that ${\eps^{-1} \in \Set{Z}}$. 

\begin{lem}
The spectral measure $\sigma_f$ of a function $f$ measurable up to $n_0$-th level partition 
is given by the Riesz product 
\begin{equation}
	\sigma_f = |\hat f|^2 \cdot \prod_{n \ge n_0} |\cP_n(t)|^2, 
\end{equation}
where
\begin{equation}
	\cP_n(t) = \frac1{\sqrt{q_n}} \sum_{y=0}^{q_n-1} e^{2\pi i\, t\omega_n(y)} 
\end{equation}
with exponential staircase frequency function 
\begin{equation}
	\omega_n(y) = m_n \frac{q_n}{\eps_n^2} (e^{\eps_n y/q_n} - 1). 
\end{equation}
\end{lem}

Remark that in $|\cP_n(t)|^2$ we can omit $-1$ in the parenthesis, and 
recall that ${h_{n+1} \ge q_n h_n}$, hence, $h_n$ are uniquely determined by~$q_n$. 
Let us define ${Q_n(t) = |P_n(t)|^2}$. We know that if we restrict $Q_n(t)$ to $(a_n,b_n)$ then 
\begin{equation}
	\left\| Q_n(t)_{(a_n,b_n)} - 1 \right\|_1 = O\left( \frac{b_n-a_n}{\sqrt{a_n m_n}} \ln(b_n m_n) \right). 
\end{equation}
Thus to apply machinery described in the beginning of this section we need estimate 
\begin{equation}
	\sum_{n=1}^\infty b_n\frac{(\ln m_n)^{1/2}}{m_n^{1/4}} < \infty, 
\end{equation}
where for simplicity we assume that ${a_n = b_n^{-1}}$. 

\begin{thm}
Let us consider the exponential staircase construction of rank one flow $T^t$ 
given by the sequence of functions
\begin{equation}
	\omega_n(y) = m_n\frac{q_n}{\eps_n^2}e^{\eps_n y/q_n}, \qquad 
	h_n = \frac{m_n}{\eps_n}. 
\end{equation}
Suppose that for an $\alpha$ with ${0 < \alpha < 1/4}$ 
\begin{equation}
	h_n^{1+\alpha} \le q_n, \qquad m_n \le h_n^{1/2-\alpha}, \qquad 
	\sum_{n=1}^\infty \frac1{m_n^{1/4-\alpha}} < \infty, 
\end{equation}
and $q_n$ has the following property: 
\begin{equation}
	\forall k \in \Set{Z} \cap 
		\left[ \frac{m_n^{1-\alpha}}{\eps_n}, \frac{m_n^{1+\alpha}}{\eps_n} e^{\eps_n} \right] 
	\quad \left\| (q_n+1) \frac1{\eps_n}k\ln k \right\|_{\Set{Z}} \le \eps_n, 
\end{equation}
where
\begin{equation}
	\| x \|_{\Set{Z}} = \min\{|\ell - x\| \where \ell \in \Set{Z}\}. 
\end{equation}
Then the flow $T^t$ has Lebesgue spectral type. 
\end{thm}

\begin{proof}
The theorem follows directly from lemmas \ref{LemSmInt} and \ref{LemConvToLOne} 
as well as ergodicity of the flow. Let us check necessary conditions 
$m_n$, $q_n$ and $h_n$ should satisfy. First, notice that it follows from
inequality ${q_n \ge h_n}$ that 
$
	h_n \ge \const \cdot 2^{c 2^n}. 
$
Since ${h_n = \frac{m_n}{\eps_n}}$, 
condition ${m_n b_n \le \eps_n^{-1}}$ follows from comparison of $m_n$ and $h_n$, 
\begin{equation}
	m_n \le h_n^{1/2-\epsilon}, \qquad 
	\frac1{\eps_n} = \frac{h_n}{m_n} \ge h_n^{1/2+\epsilon}, 
\end{equation}
and ${\eps^{-1} \le h_n \le q_n^{1/(1+\eps)}}$. 
Thus, the series $\sum_{n=1}^\infty b_n\frac{(\ln m_n)^{1/2}}{m_n^{1/4}}$ converges 
and all other requirements of $m_n$, $\eps_n$ and $q_n$ are satisfied, so 
we can apply the above flatness lemma \ref{LemSmInt} concerning 
exponential staircase trigonometric sums 
with convergence lemmas \ref{LemConvOfTrigSums} and \ref{LemConvToLOne}. 
To complete the proof we have to mention that it follows from ergodicity 
of our rank one flow that the spectral measure has no atom at zero. 
\end{proof}